\documentclass[10pt,a4paper]{amsart}
\usepackage[english]{babel}
\usepackage[utf8x]{inputenc}
\usepackage{ucs}
\usepackage{tikz}
\usepackage{amsmath}
\usepackage{amsthm}
\usepackage{amsfonts}
\usepackage{amssymb}
\usepackage{amsthm}
\usepackage{hyperref}
\usepackage{dsfont}
\usepackage{mathrsfs}
\usepackage{eqnarray}
\write18{}
\usepackage{xcolor}
\usepackage{graphicx}
\usepackage{pdfpages} 
\usepackage[nothm]{thmbox}
\usepackage[pagewise]{lineno}

\hypersetup{
    colorlinks=true,
    linkcolor=red,
    filecolor=magenta,      
    urlcolor=cyan,
}
\numberwithin{equation}{section}

\newcommand{\inte}[1]{\underset{#1}{\int}}
\newcommand{\somme}[1]{\underset{#1}{\sum}}
\newcommand{\limi}[1]{\underset{#1}{\lim}}
\newcommand{\tends}[1]{\underset{#1}{\longrightarrow}}

\newcommand{\fonction}[5]{\begin{array}{r r c l}
					#1 \hspace{2mm} : & #2 & \to & #3 \\
					& #4 & \mapsto & #5 \\
			   \end{array}}
\newcommand{\fonctionbis}[4]{\begin{array}{r c l}
					   #1 & \to & #2 \\
					 #3 & \mapsto & #4 \\
			   \end{array}}	  
\newcommand{\quotient}[2]{{\raisebox{.2em}{$#1$}\left/\raisebox{-.2em}{$#2$}\right.}}

\newcommand{\resp}[1]{\textit{(resp. }#1\textit{)}}

\DeclareMathOperator{\vol}{vol}
\DeclareMathOperator{\lk}{lk}
\DeclareMathOperator{\spn}{vect}

\theoremstyle{definition}

\newtheorem{definition}[equation]{Definition}
\newtheorem{remark}[equation]{Remark}

\theoremstyle{theorem}

\newtheorem{theoreme}[equation]{Theorem}
\newtheorem{proposition}[equation]{Proposition}
\newtheorem{corollary}[equation]{Corollary}

\newtheorem{lemma}[equation]{Lemma}

\title{A spectral approach to the linking number in the 3-torus}

\author{Adrien Boulanger}

\date{}

\begin{document}

\maketitle

\begin{abstract}
Given a closed Riemannian manifold  and a pair of multi-curves in it, we give a formula relating the linking number of the later to the spectral theory of the Laplace operator acting on differential one forms. As an application, we compute the linking number of any two multi-geodesics of the flat torus of dimension 3, generalising a result of P. Dehornoy.
\end{abstract}

\section{Introduction.}

Let $(M,g)$ a closed Riemannian manifold of dimension 3. We call \textbf{a curve} an embedding of the oriented circle. A \textbf{multi-curve} is a finite collection of disjoint curves.
We say that a multi-curve is \textbf{homologically trivial} if its homology class vanishes, as a cycle of $M$. \\

Recall that given two homologically trivial multi-curves $\Gamma, \Upsilon$, one can define their \textbf{linking number} in taking any surface $S_{\Gamma}$ whose boundary is $\Gamma$ and (algebraically) in intersecting it against $\Upsilon$ ;
		$$ \lk( \Gamma, \Upsilon ) := i( S_{\Gamma}, \Upsilon)  \ .$$

	\begin{figure}[!h]
	\begin{center}
			\def\svgwidth{0.35 \columnwidth}
			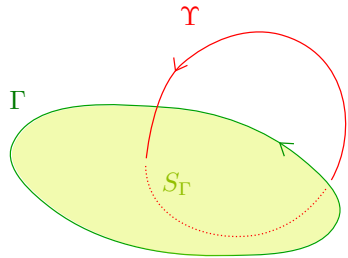
	\caption{ Here, both collections $\Gamma$ and $\Upsilon$ consist of a single curve. Their linking number is $\pm 1$, depending on the global orientation.}
	\end{center}
	\end{figure}

It is not immediate that this number is well defined, because of the choice involved about a surface $S_{\Gamma}$. As a general reference to the notion of linking number, one can recommend \cite[Chapter III Section 4]{livrearnold} and \cite[Section 28]{livrebott}. Our main theorem relates the linking number with the spectral theory.

\begin{theoreme}
	\label{theo enlacement spectrale}
	 Let $(M,g)$ be a closed Riemannian manifold and $\Gamma$, $\Upsilon$ two disjoint homologically trivial multi-curves, they link according to the following formula:
	\begin{equation}
	\label{theo enlacement spectrale equation}
			\lk(\Gamma, \Upsilon) = \limi{t \to 0} \ \somme{k \ge 0} \ e^{-\lambda_k t} \inte{\Gamma} \eta_k \inte{\Upsilon} *\left( \frac{d \eta_k}{\lambda_k} \right) \ ,
		\end{equation}
	where $(\eta_k)_{ k \in \mathbb{N}}$ denotes an eigenvector basis with corresponding eigenvalues $(\lambda_k)_{k \in \mathbb{N}}$ of the Laplace operator $\Delta$ acting on the Hilbert space of square integrable 1-differential forms in $\ker(\Delta)^{\perp}$.
\end{theoreme} 

Note that this theorem relates a topological  number with metric quantities. In particular the right member of Formula \ref{theo enlacement spectrale equation} does not depend on the underlined metric $g$. \\

Theorem \ref{theo enlacement spectrale} can be used if one has enough knowledge of the spectral theory of $(M,g)$, as it is the case of the canonical flat torus $\mathbb{T}^3$. We prove a general formula for the linking number of multi-curves consisting of geodesics of $\mathbb{T}^3$. However, not to burden this introduction, we postpone the statement in Section \ref{3}. Specialising our formula to the case of closed orbits of the geodesic flow on the 2-torus $\mathbb{T}^2$ gives the 

\begin{corollary}
\label{Pierre2}
Let $\Gamma = (\gamma^i)_{I \in I}$ and $\Upsilon = (\upsilon^j)_{j\in J}$ two homologically trivial multi-curves in  $\mathbb{T}^3$ consisting of periodic orbits of the $\mathbb{T}^2$ geodesic flow. They link according to the formula:
$$ \lk (\Gamma, \Upsilon) = \somme{i \in I,j \in J}   i ( \gamma^i, \upsilon^j ) \frac{1 -  \frac{x_{i,j}}{\pi}}{2} \ ,$$
where $x_{i,j}$ denotes the unique determination in $[0,2\pi[$ of the oriented angle $\theta$ made at each intersection point (see Figure 2) and $i ( \gamma^i, \upsilon^j )$ denotes the algebraic intersections between the projections on $\mathbb{T}^2$ of the curves $ \gamma^i$ and $\upsilon^j$.
\end{corollary}

\begin{figure}[h!]
\begin{center}
	\def\svgwidth{0.3 \columnwidth}
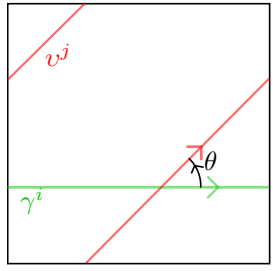
\caption{Here, the intersection number is 1. The angle $\theta$ defined in Corollary \ref{Pierre2} is represented in black.}
	\end{center}
	\label{fig tore pierre}
\end{figure}

Another formula was found by P. Dehornoy using different methods in his PhD \cite{thesedehornoy}. Our formula shows, in a clear way, that the linking number entertains some interactions with the intersection number on the curves projected on the basis.\\ 

We now briefly survey old and more recent statements about the linking number. The first occurence of the notion of linking number goes back to Gauss's studies on electromagnetism (see \cite{articlericcahistoire}). Gauss noticed that integrating the magnetic field generated by an electric power flowing in a close wire $\gamma$ - for us a differential form $\omega_{\gamma}$ - along any closed curve $\upsilon$ gives a number which does not depend on the homology class of $\upsilon$ in the complementary of $\gamma$. That is to say, the differential form $\omega_{\gamma}$ is closed. \newline 

In fact, Gauss went further in his study:  he gives in $\mathbb{R}^3$ an explicit formula expressing the differential form $\omega_{\gamma}$. Let $x \notin \gamma$ and $X(x) \in T_x(\mathbb{R}^3)$, then:
\begin{equation}
	\label{equationprimite gamma}
		(\omega_{\gamma})_x(X) = \frac{1}{4\pi} \inte{[0,2\pi]} \det \left( \gamma'(s), X(x) , \frac{\gamma(s) - x}{||\gamma(s) - x||^3} \right) ds  \ ,
\end{equation}

where $\gamma(s)$ denotes any parametrisation compatible with the curve $\gamma$ orientation, and $ || \cdot ||$ is the Euclidean norm of $\mathbb{R}^3$. Back to this days, there was no topological definition of the linking number so that, following Gauss, one could have defined it setting:
\begin{equation}
	\label{formuleGausslk}
	 \lk(\gamma,\upsilon) := \inte{ \upsilon} {\omega_{\gamma}} = \frac{1}{4\pi} \inte{[0,2\pi]} \inte{[0,2\pi]} \det \left( \gamma'(s), \upsilon'(t) , \frac{\gamma(s) - \upsilon(t)}{||\gamma(s) - \upsilon(t)||^3} \right) ds \ dt \ . 
\end{equation}

Gauss's formula had been related later on with the linking number defined at first, see for example \cite[Chapter 3, Section 4]{livrearnold}. It is still an active research field to try to get Gauss-like formulas for the linking number and its natural generalisation \cite{articledeturck3noeuds} \cite{articledeturcknsphere}. \\

Formula \eqref{formuleGausslk} also suggests the existence of a universal object which, integrated over a pair of homologically trivial multi-curves, gives back their linking number. A \textbf{linking form} $\Omega$ is an integrable $(1,1)$-differential form satisfying for any two homologically trivial disjoint multi-curves $\Gamma$ and $\Upsilon$ 
	$$ \lk(\Gamma, \Upsilon) = \inte{\Gamma} \inte{\Upsilon} \Omega \ .$$ 

One can think of a $(1,1)$-differential form as a 2-differential form; we will get back on the $(1,1)$-form precise definition in Subsection \ref{2.2}. \\

The definition of linking form was introduced by Arnold (see \cite[Chapter III, Section 4]{livrearnold}) to generalize Moffatt's \cite{articlemoffatt} interpretation of the helicity. Let us recall briefly how is defined the later. \\

Let $X$ be a vector field preserving a probability measure $\mu$ in the Lebesgue class whose asymptotic cycle vanishes. This assumption implies that the 2-differential form $i_X \mu$ is exact, meaning that there is a 1-differential form $\alpha$ such that $ d \alpha = i_X \mu$. On can show that
	$$ \mathcal{H}(X) := \inte{M} \alpha \wedge d \alpha $$
 does not depend on the choice involving the primitive $\alpha$. We call this number the helicity of the vector field $X$. This notion was introduced by \cite{articlemoreau} and \cite{articlewoltjer} to study certain energies associated to vector fields solution of some partial differential equations. Note that the asymptotic cycle assumption is automatically satisfied in some natural situations, for example when the ambient manifold is a homology sphere or if $X$ is the Reeb flow associated to a contact structure.\\

Arnold interpreted the helicity of a vector field $X$ as some average of the asymptotic linking number of two trajectories of the flow. Given $x,y \in M$, we consider the trajectories starting off $x$ and $y$ of the flow $X$ at times $t$ and $s$. We close them by gluing their extremities using a small path, we compute the linking number, we divide by the product $t s$ and one would like to let $s,t \to \infty$. To do it, one needs to show this limit to be almost everywhere well defined; this is one of the reason why Arnold introduced the notion of linking form. Actually, he showed that the linking form is \textbf{integrable}, see Remark \ref{rem parametrization}, which allows one to define the previous limit using Birkhoff's ergodic Theorem. Arnold recovered the helicity from the construction of the linking form itself; as the kernel of some inverse of the differential operator $d$, acting on 1-differential forms. See \cite[Chapter III, Section 4]{livrearnold} for more detail. This perspective on the helicity had been developped in \cite{articlevogellinking}, \cite{articlevogelcontact} and \cite{articledeturckgausshel}. \\

Arnold also noticed that linking forms always exists on compact manifolds. This had been precised by T.Vogel in \cite{articlevogellinking}, relying on G. de Rham's work on Hodge theory \cite[Section 28]{livrederham}. We denote by $g^1(x,y)$ the kernel of the \textbf{Green operator}, the inverse operator of the Laplace one. We have the

\begin{theoreme}{\cite[Theorem 3]{articlevogellinking}} 
\label{theo vogel de rham}
Let $(M,g)$ a compact Riemannian manifold. The $(1,1)$-differential form 
	\begin{equation}
		\label{eqtheo vogel de rham}
		 \Omega(x,y) = *_y d_y g^1(x,y) 
	\end{equation}
is an integrable linking form. We call this linking form the \textbf{de Rham-Vogel's linking form}.
\end{theoreme}

T. Vogel's proof relies on Arnold's remark that any inverse operator of $d$ gives rise to a linking form, up to Hodge duality. This theorem shows the existence of linking forms on closed manifolds, but does not come with a simple formula like Gauss's one \eqref{formuleGausslk}. There is, up to the author's knowledge, only two others known formulas of this type, found in \cite{articledeturck}. The first one holds for the hyperbolic 3-space and the second one for the round 3-sphere. The authors find such a formula in exhibiting a "fundamental solution of Maxwell's Equations", meaning in exhibiting the de Rham-Vogel linking form defined above. \\

\textbf{Outline of the paper.} Formula \eqref{eqtheo vogel de rham} relates the linking number to the Green operator associated to the Riemannian structure on the underlined manifold $M$. Thanks to the \textbf{spectral theory}, providing $M$ to be closed, one knows that the Green operator is diagonalisable. It is therefore natural to try to get a spectral-linking formula in splitting the Green operator with respect to a basis of 1-eingenforms of the Laplace operator. \\

However, the integration current over a closed curve being not regular enough - not square integrable - one can not hope to readily get this spectral-linking formula. To circumvent this difficulty, in order to reach some more regularity, we will use the \textbf{heat operator} to diffuse the integration currents. This is the heart of Section \ref{3} which lasts with the proof of our Theorem \ref{theo enlacement spectrale}. \\

One wants then to use Theorem \ref{theo enlacement spectrale} to compute the linking number of some collections of curves in manifolds for which the spectral theory is well known, like the canonical flat 3-torus $ \mathbb{T}^3 := \quotient{\mathbb{R}^3}{\mathbb{Z}^3}$. \\

We will give in Section \ref{4} a general formula to compute the linking number between two homologically trivial collections of multi-curves consisting of geodesics of $\mathbb{T}^3$. Natural collections of geodesics in $\mathbb{T}^3$ are given by closed orbits of the geodesic flow of $\mathbb{T}^2 := \quotient{\mathbb{R}^2}{\mathbb{Z}^2}$, its unitary tangent bundle being identified with $\mathbb{T}^3$, see Section \ref{4.4}. Specialising our formula in this case, we find a new formula which computes the linking number between two homologically trivial collections of the $\mathbb{T}^2$-geodesic periodic orbits. Our formula, as well as P. Dehornoy's one, shows that the linking number of two homologically trivial multi-curves consisting of periodic orbits of the $\mathbb{T}^2$-geodesic flow have always the same sign. This question of the sign is closely related to the notion of Birkhoff sections. For more detail about this recent aspect on the linking number, see Subsection \ref{4.5} and references therein.\\

\textbf{Acknowledgement.} I would like to thank P. Dehornoy for the many interesting conversations we had around this article and more. I also want to thank my PhD advisors, G. Courtois and F. Dal'Bo, for pointing out the question as well as for their comments about the redaction of this article. 

\section{Green operators and linking forms}

\label{2}

This section is devoted to introduce all the objects we will use later on.

\subsection{The Green kernel}
\label{2.1}

Let $(M,g)$ be a closed manifold of dimension $p$. We denote by  \\
\begin{itemize}
\item $\mu_g$, the volume form associated to the metric $g$;
\item $\Omega^*(M) = \underset{ 0 \le k \le p}{\oplus} \Omega^k(M)$, the space of all differential forms, split with respect to the degree $k$;
\item $*$, the Hodge operator, or Hodge star, which satisfies the following identity
	\begin{equation}
		\label{equ hodge star}
			** = (-1)^{k(p-k)} \ ;
	\end{equation}
\end{itemize}
\textit{Note that we abuse the notation in omitting the degree $k$ of the underlined differential form.}
\begin{itemize}
\item $d$, the exterior differential operator on $\Omega^*(M)$. \\
\end{itemize}   

The Hodge star is defined to endow the vectorial space $\Omega^k(M)$ with a scalar product:
	\[ \big<\alpha \cdot \beta \big> =  \inte{M} \alpha \wedge * \beta  \ . \]
	
With respect to it, the operator $d$ has a unique adjoin operator, denoted by $\delta$, satisfying by definition:
    \[ \big< d \alpha \cdot \beta \big> = \big< \alpha \cdot \delta \beta \big> \ .\]
    
A straighforward computation involving the Hodge star definition and the Stokes formula gives 
\begin{equation}
	\label{equ delta}
	\delta = (-1)^{p(k +1) +1 } *d* \ . 
\end{equation}

We now have all the material required to define the Laplace operator.

\begin{definition}
The \textbf{Laplace operator}, denoted by $\Delta \curvearrowright \Omega^*(M)$, is defined by
	$$\Delta := d \delta + \delta d \ . $$
\end{definition}

Note that the Laplace operator stabilises all differential forms spaces of fixed degree, and that it is self-adjoin with respect to the scalar product $ \big< \cdot, \cdot \big>$. We denote by $\Delta^k$ its restriction to $\Omega^k(M)$. \newline 

A differential form $\alpha$ is said to be harmonic if it lies in the kernel of the Laplace operator, denoted by $\ker \Delta$. The space of all harmonic $k$-forms being identified, from a famous theorem of Hodge, to the $k$-nth homology group of $M$ (see \cite[page 46]{livrerosenberg} for example), the Laplace Operator $\Delta^k$ can not be invertible in general. However, nothing prevents it to be invertible in restriction to the orthogonal space of its kernel. We denote by $\pi_k$ the orthogonal projection on $\ker( \Delta^k )^{\perp}$.

\begin{definition}
	\label{defnoygreen}
	A \textbf{Green operator}, denoted by $G^k$, is any operator satisfying the following equation on the space of smooth differential forms of degree $k$:
		\begin{equation}
			\label{eqdef noyau de Green} 
				G^k \circ \Delta^k = \Delta^k \circ G^k = Id - \pi_{\mathcal{H}^k}
		\end{equation}				
\end{definition}

Such an operator always exists, providing $M$ to be closed (see for example \cite[Section 3]{livrederham}). There is a slight ambiguity about $G^k$, which is fully determined up to its restriction on the space $\ker( \Delta^k )$. From now, we will suppose that $G( \ker( \Delta^k ) ) = \{ 0 \}$, allowing one to speak of \textbf{the} Green operator. \\

Green operators are \textbf{kernel operators}, meaning that there is a smooth family of endomorphism - what we call a (1,1)-form - $g^k(x,y) : \Lambda^k(T_x M)  \to \Lambda^k(T_y M)$, indexed by $M \times M \setminus \text{Diag}$ such that for all smooth differential forms $\alpha$ of degree $k$ 

\begin{equation*}
		 G^k(\alpha)_y = \inte{y  \in M} g^k(x,y)(\alpha_x) d \mu_g(x) \ .
\end{equation*}

\subsection{Differential (1,1)-forms} 
\label{2.2}
We give in this subsection the precise definition of a (1,1)-form. We also explain how to integrate them over a pair of multi-curves. Given an euclidean space $E$, we denote by $\sharp$ the musical endomorphism which maps some vector $X \in E$ on its dual linear form, so in $E^*$, according to the Euclidean structure on $E$. 

\begin{definition}
	Let $M$ be a manifold, We call a \textbf{(1,1)-form} a family of morphisms $T^*_x(M) \to T^*_y(M)$ indexed by $M \times M$.
\end{definition}

 Let $\gamma$ and $\upsilon$ two parametrised curves by $s$ and $t$. We define the integral over a pair of curves of a (1,1)-form $\Omega$ as 
	$$ \inte{\gamma} \inte{\upsilon} \ \Omega \ :=  \inte{\gamma} \left( \inte{0}^1 \Omega( \upsilon(s) , y ) \left( \upsilon'(s)^{\sharp} \right) \ ds \right) \ . $$

Moreover, the following integral - an element of $T^*_y(M)$ -
$$ \inte{0}^1 \Omega( \upsilon(s), y) \left( \upsilon'(s)^{\sharp} \right) \ ds  $$

does not depend on a choice of parametrisation, since $\Omega(x,y)$ is linear. So that we will prefer to denote it for short as
	$$ \inte{\upsilon} \Omega \left((\cdot)^{\sharp},y \right) \ ,$$
omitting the underlined parametrisation. \\

This formula shows clearly that the linking form enjoys some bilinearity. In fact, if we denote $\Upsilon = \underset{i \in I}{\cup} \upsilon^i $ and $\Gamma = \underset{j \in J}{\cup} \gamma^j $ we have:
\begin{align*}
		\lk( \Gamma, \Upsilon) & = \inte{\Gamma \times \Upsilon} \Omega 
 		= \inte{\cup \gamma_i \times \cup \gamma_j} \Omega 
 		= \somme{i \in I ,  j \in J} \inte{ \gamma_i \times \gamma_j} \Omega \ .
\end{align*}

Note that we did not ask either the curve $\gamma_i$ or the curve $\upsilon_j$ to be homologically trivial. \\

\subsection{The de Rham-Vogel's linking form}

\label{2.3}

Recall that the De Rham-Vogel's linking form is defined as 
	$$ *_y d_y g^1(x,y) \ ,$$ 
which may be slightly confusing at first. What does it mean to consider the image by $*_y d_y$ of a family of morphism from $T^*_x M $ to $T^*_y M$ ? \\

Given $\alpha \in T_x^*(M)$, the Green kernel defines a differential form by:
	$$ \alpha_y := y \mapsto g(x,y)\big( \alpha(x) \big) \in T_y^*(M) \ .$$

This differential form is smooth on $M \setminus \{x\}$, which allows one to take its image by the operator $*d$ wherever it makes sense. This gives rise to an other linear morphism 
	$$ \fonctionbis{T_x^*(M)}{T_y^*(M)}{\alpha}{\alpha_y \ ,} $$
which turns out to correspond to the kernel of the operator $ \alpha \mapsto (* d G^1)(\alpha)$. Then, in the end, De Rham's notation $*_y d_y g^1(x,y)$ is to be understood as the kernel of the operator $ \alpha \mapsto (* d G^1)(\alpha) \ .$
	
\begin{remark} \
	\begin{itemize}
	\item As pointed out by Arnold, any kernel associated to the inverse operator of the exterior differential $d$ is a linking form. The operator $\alpha \mapsto (* d G)(\alpha) $ is actually one of them, up to Hodge duality. See \cite[Lemma 2]{articlevogellinking}.
	\item The singularity of the (1,1)-form $g^1(x,y)$ along the diagonal is roughly equivalent to $r^{-1}$. Thus, after one differentiation, this singularity turns to be in $r^{-2}$, which is still integrable in dimension three, see \cite[Theorem 23 page 134]{livrederham}. So that what we meant by integrable is that for every $x$ the function 
	$$ y \mapsto || *_y d_y g(x,y) || $$
is integrable on $M$ with respect to $\mu_g$. The notation $|| \cdot ||$ stands for the linear morphism norm induced by the metric $g$. 
	\end{itemize}
\end{remark}
	
\subsection{Behaviour of the linking form under isometries.}
\label{2.4}

The de Rham-Vogel linking form being constructed from a metric, it is natural to look into how it behaves under an isometry $\Phi$. The isometry  $\Phi$ commutes with the Hodge star as well as with the exterior differential $d$. In particular, it commutes with every operators made out of this two ones, as the Laplace operator and its inverse, the Green operator. Looking at the kernel of the later, this commutation relation can be read as:
	$$ (\Phi_1)_* g^1( x,y) = (\Phi_2)^* g^1( x, y) \ ,$$ 

where $\Phi_1$ \resp{$\Phi_2$} denotes the $\Phi$-action on the first \resp{second} factor of the product $M \times M$. In particular, the diagonal action of $\Phi$ on the product $M \times M \setminus \text{Diag}$ preserves the Green kernel, and thus the de Rham-Vogel's linking form. Since we will use this remark to simplify a bit the calculations performed in Section \ref{4}, we set it in a 

\begin{proposition}
\label{translation}
Let $\gamma$ and $\upsilon$ two curves, not necessary homologically trivial, $\Phi$ an isometry of $(M,g)$, then :
$$ \inte{\gamma} \inte{\upsilon} *_y d_y g^1 = \inte{\Phi^{-1}(\gamma)} \inte{\Phi^{-1}(\upsilon)} *_y d_y g^1 \ . $$
\end{proposition}

\section{The spectral-linking formula}
\label{sec an other formula}
\label{3}

Let us recall that $\Delta$ is self-adjoin with respect to $\big<\cdot, \cdot\big>$. It is well known that self-adjoin operators are diagonalisable in finite dimension; it is actually still the case for the Laplace operator, providing the underlined manifold to be closed. 

\begin{theoreme}\cite[Theorem 1.30]{livrerosenberg}
Let  $(M,g)$ be a closed Riemannian manifold. There is a orthonormal basis $(\eta_n)_{n \in \mathbb{N}}$ of 1-differential forms, meaning that $ \big< \eta_i,\eta_j \big>  = \delta_i(j)$, and a sequence of non negative numbers $(\lambda_n)_{n \in \mathbb{N}}$ such that
			$$ \Delta \eta_n = \lambda_n \eta_n \ .$$
\end{theoreme}

In particular, if $\alpha \in \ker (\Delta)^{\perp}$ we have:
		$$ \alpha = \somme{n \in \mathbb{N}} \big<\eta_n \cdot \alpha \big> \eta_n \ . $$

Formally, one would like to write the Green operator as 
	$$ g^1(x,y) := \somme{n \in \mathbb{N}} \frac{1}{\lambda_n} \eta_n(x) \otimes \eta_n(y) \ , $$
	
which gives the following expression for the de Rham-Vogel's linking form
	$$ *_y d_y g^1(x,y) := \somme{ n \in \mathbb{N}} \ \eta_n(x) \otimes * d \left( \frac{ \eta_n(y)}{\lambda_n} \right) \ . $$

Keeping it formal, one would like then to integrate each factor along $\gamma$ and $\upsilon$ to get:
	$$	\lk(\gamma, \upsilon) =   \ \somme{n \in \mathbb{N}} \  \inte{\gamma} \eta_n \inte{\upsilon} *\left( \frac{d \eta_n}{\lambda_n} \right)  \ .$$ 

However, the previous series does not converge \textit{a priori}. In fact, an integration current over a curve is not square integrable and therefore cannot be decomposed with respect to the orthonormal basis $(\eta_n)$. To circumvent this difficulty, we will regularise them thanks to the use of the heat kernel, from which the term $e^{- \lambda_n t}$ of formula \ref{theo enlacement spectrale equation} comes from. As a corollary of this approach, we are able to prove the following stronger version of Theorem \ref{theo enlacement spectrale}.

\begin{theoreme}
	\label{theo formule localisebis}
	 Let $(M,g)$ be a closed Riemannian manifold and  $\Omega$ the de Rham-Vogel's linking form, then for all pair of curves $\gamma$ and $\upsilon$, not necessary homologically trivial we have:
	\begin{equation}
	\label{equation localisebis}
			\inte{\gamma} \inte{\upsilon} \Omega = \limi{t \to 0} \ \somme{k > 0} \ e^{-\lambda_k t} \inte{\gamma} \eta_k \inte{\upsilon} *\left( \frac{d \eta_k}{\lambda_k} \right) \ . 
	\end{equation}
where $(\eta_k)_{ k \in \mathbb{N}}$ denotes an eigenvector basis of the Laplace operator $\Delta$ acting on the Hilbert space $\ker(\Delta)^{\perp}$ - viewed as a subspace of square integrable differential forms - and $(\lambda_k)$ the associated eigenvalues.
\end{theoreme}

All the rest of this section is dedicated to the proof of the above theorem.

\subsection{The heat operator on 1-differential forms}

\label{3.1}

The following definition is the key to regularise the integration currents. More detail about generalised heat kernels can be find in \cite[Section 2.3]{livrevergne}.

\begin{definition}
	\label{definition heat kernel}
Let $(M,g)$ be a closed Riemannian manifold and $\eta$ be a continuous, bounded 1-differential form. The following Cauchy problem of unknown $(\eta_t)_{t \in\mathbb{R_+}}$ 
	\begin{equation*}
		\left\{
		\begin{split} 
				 \Delta \eta_t & +  \partial_t \eta_t = 0\\
				 \eta_0 & = \eta \\
		\end{split}
		\right.
	\end{equation*}
has a unique solution. We denote by $e^{-t\Delta^1}$ the \textbf{heat operator} which maps $\eta$ to the time $t$ solution of the above Cauchy problem.We denote by $p^1_t$ the \textbf{heat kernel}, which satisfies by definition:
	$$ (\eta_t)_y = \inte{M} \  p^1_t(x,y)(\eta_x) \ d \mu_g(x) \ . $$
Moreover, one has 
	 $$ e^{-t\Delta^1}(\eta) \tends{ t \to 0} \eta $$ 
for the uniform convergence topology. 
\end{definition}

In particular, if $U$ and $V$ are two closed disjoin subsets of $M$, one has 
$$ p^1_t(x,y) \tends{ t \to 0} 0 $$
uniformly on $U \times V$. \newline

The heat kernel has the interesting property to be smooth for all $t > 0$, on the contrary of the Green operator. In particular, it is decomposable according to an orthonormal basis of eigenforms. 

\subsection{The diffused curves}
\label{3.2}

Let $\gamma$ be a curve of $M$. We denote by  $L^1(\Omega^1(M))$ the space of integrable 1-differential form, meaning forms whose punctual norm is integrable over $M$ with respect to the Riemannian measure.

\begin{definition}
We call the $\gamma$-\textbf{diffused curve}, denoted by $e^{-t\Delta^1}(\gamma)$, the following family of linear forms indexed by $t > 0$:
	\begin{equation*}
		\fonction{e^{-t\Delta^1}(\gamma)}{L^1(\Omega^1(M))}{\mathbb{R}}{\beta}{ \inte{\gamma} e^{-t\Delta^1}(\beta) \ . } 
	\end{equation*}	 
\end{definition}
	
This diffusing process associates to each $t>0$ a differential form approximating the integration current over the curve $\gamma$: the smaller is $t$, the better is the approximation.

\begin{lemma}
	\label{lemma diffused curve}
For all  $\beta \in L^1(\Omega^1(M))$ continuous on a neighbourhood of $U$ of the curve $\gamma$, we have 
		\begin{equation*}
			e^{-t\Delta^1}(\gamma)(\beta) \tends{ t \to 0} \inte{\gamma} \beta \ .
		\end{equation*}
\end{lemma}

\textbf{Proof:}  we have been careful to consider a differential form $\beta$ integrable. So, since the heat kernel converges uniformly to $0$ away from the diagonal, we have:
	$$ \big| e^{-t\Delta^1}(\gamma)(\beta) - \inte{\gamma}  \inte{U} \ p_t^1(x,y) (\beta_x) \ d\mu_g(x) \big| \tends{t \to 0} 0 \ .$$
	
The differential form $\beta$ being continuous on $U$, from the very definition on the heat kernel we have: 
	$$  \inte{U} \ p_t^1(x,y) (\beta_x) \ d \mu_g (x) \tends{t \to 0} \beta_y $$

uniformly. Therefore, one is allowed to permute limit and integral to get 
	$$ \inte{\gamma}  \inte{U} p_t^1(x,y) ( \beta_x) \ d \mu_g (x) \tends{t \to 0} \inte{\gamma} \beta \ , $$

which is the expected result. \hfill $\blacksquare$ \newline 

If now $\upsilon$ is curve disjoin of $\gamma$, recall that the 1-differential form 
	$$ {(\omega_{\upsilon})}_y := \inte{\upsilon} \Omega \big((\cdot)^{\sharp}, y \big) $$
	
is integrable, where $\Omega =*_1 d_1 g_1$ is the de Rham-Vogel's linking form. Applying Lemma \ref{lemma diffused curve} gives readily the

\begin{corollary}
	\label{corollary lemma diffused}
	For any two curves $\gamma$ and $\upsilon$ we have  
		$$ e^{-t\Delta^1}(\gamma)(\omega_{\upsilon}) \tends{t \to 0} \inte{\gamma} \inte{\upsilon} \Omega \ .$$ 
\end{corollary}

The goal is know to identify, $t > 0$ being fixed, the left member of the above equation to the series appearing in Equation \eqref{equation localisebis}. We will conclude by summoning the above corollary to recover Theorem \ref{theo formule localisebis} by letting $t \to 0$.

\subsection{The approximating series} \label{subsec approximate series} 

\label{3.3} The benefits of having diffused the integration current is to allow one to write the left member of \eqref{equation localisebis} as a scalar product of two smooth 1-differential forms. We will conclude in using Plancherel's formula, allowing one to write down with respect to an orthonormal basis this scalar product. \\

\begin{lemma}
	\label{lemme PS form noyau}
	For all differential forms $\beta \in L^1(\Omega^1(M))$ and all $t > 0$ we have 
		$$ e^{-t\Delta^1}(\gamma)(\beta) = \left< \beta \cdot  \inte{\gamma} p^1_t \big( (\cdot)^{\sharp} , y \big)\right> \ .  $$
\end{lemma}

Note that the scalar product is well defined since the differential form $\inte{\gamma} p^1_t \big( (\cdot)^{\sharp} , y \big)$ is smooth. \\

\textbf{Proof:} the operator $e^{-t \Delta}$ being self-adjoin and since $i_X(\alpha)(x) = g_x \left( X^{\sharp} \cdot \alpha \right)$, we have the following identity for any 1-differential forms $\beta$ and any vector field $X$: 
	$$ i_X(y) \big( p_t(x,y) \beta_x \big)  = g \left( \beta_y \cdot \left( p_t(x,y) (X^{\sharp}_x )\right) \right) \ .  $$

Therefore, setting $X_{x} = \gamma'(s)$ and integrating along $\gamma$, one gets:
$$ \inte{\gamma} p_t^1(x, \cdot)(\beta_x) = g_y \left( \beta_y, \inte{\gamma} p^1_t \big( (\cdot)^{\sharp} , y \big) \right) \ ,$$

which gives, after integration over $M$ with respect to $\mu_g$,
	$$ \inte{M} \inte{\gamma} \ p_t^1(x, y)(\beta_x) \ d \mu_g(y) = \left< \beta \cdot  \inte{\gamma} p^1_t \big( (\cdot)^{\sharp} , y \big) \right> \ . $$

We conclude recalling that the form $\beta$ is integrable, which allows one to switch both integrals of the above equation left member, recovering our definition of a diffused curve. \hfill $\blacksquare$ \\

We conclude the proof of \ref{theo formule localisebis} as announced by identifying the right member of \eqref{equation localisebis} with some series. 

\begin{lemma}
For all $t > 0$ we have:
	\begin{equation*}
		 	e^{-t\Delta^1}(\gamma)(\omega_{\upsilon}) = \somme{k > 0} e^{-\lambda_k t} \inte{\gamma} \eta_k \inte{\upsilon} *\left( \frac{d \eta_k}{\lambda_k} \right)  \ .
	\end{equation*}	 
\end{lemma}

\textbf{Proof:} we start with using the semi-group property of the heat operator $e^{-t\Delta}$,
	\begin{align*}
		e^{-t\Delta^1}(\gamma)(\omega_{\upsilon}) & = e^{-\frac{t}{2} \Delta^1}(\gamma) \left( e^{-\frac{t}{2} \Delta^1}(\omega_{\upsilon}) \right) \ ,
	\end{align*}

for which we apply Lemma \ref{lemme PS form noyau} to get 
	$$ e^{-\frac{t}{2} \Delta^1}(\gamma) \left( e^{-\frac{t}{2} \Delta^1}(\omega_{\upsilon}) \right) =   \left< \inte{\gamma}  p_{\frac{t}{2}}^1( (\cdot)^{\sharp},y ) \ \cdot \  \left( e^{-\frac{t}{2}\Delta^1}  (\omega_{\upsilon}) \right) \right> \ . $$
 
Both 1-differential forms appearing in the above equation being smooth, one is able to write down this scalar product with respect to an orthonormal basis consisting of the Laplace operator eingenforms: 
	\begin{equation*}
		 e^{-t\Delta^1}(\gamma)(\omega_{\upsilon}) = \somme{ k \in \mathbb{N}}  \left<  \left[ \inte{\gamma} p_{\frac{t}{2}}^1(( \cdot)^{\sharp},y ) \right]  \cdot \eta_k \right> \left< e^{-\frac{t}{2}\Delta^1}  (\omega_{\upsilon}) \cdot \eta_k \right> 
	\end{equation*}

It remains then to prove both the two following identities:
\begin{equation}
	\label{equation 1}
e^{-\frac{\lambda_k t}{2}} \inte{\gamma} \eta_k =
				 	 \left<  \left[ \inte{\gamma} p_{\frac{t}{2}}^1(( \cdot)^{\sharp},y ) \right]  \cdot \eta_k \right>
\end{equation}

\begin{equation}
	\label{equation 2}
	 	\frac{e^{-\frac{\lambda_k t}{2}}}{\lambda_k} \inte{\upsilon} *d \eta_k = \left< e^{-\frac{t}{2}\Delta^1}  (\omega_{\upsilon}) \cdot \eta_k \right> 
\end{equation}

We start off the right member of \eqref{equation 1}. Recalling Lemma \ref{lemme PS form noyau} the other way around, one gets 
 $$ \left<  \left[ \inte{\gamma} p_{\frac{t}{2}}^1(( \cdot)^{\sharp},y ) \right]  \cdot \eta_k \right> = \inte{\gamma}  e^{-\frac{t}{2}\Delta^1}\left(\eta_k \right) \ .$$
 
Then, $\eta_k$ being an eigenform of eigenvalue $\lambda_k$, we have
	\begin{equation*}
			e^{-\frac{t}{2}\Delta^1}\left(\eta_k \right) \  = e^{- \frac{t}{2} \lambda_k}  \inte{\gamma} \eta_k  \ .
	\end{equation*}
which proves that \eqref{equation 1} holds. \\

Let us show the same way that \eqref{equation 2} occurs as well. We start again from the right member:
	$$ \left< e^{-\frac{t}{2} \Delta^1}(\omega_{\upsilon}) \cdot \eta_k \right> \ . $$

The differential 1-form $\omega_{\upsilon}$ being integrable and the operator $e^{-\frac{t}{2} \Delta}$ being self-adjoin we have:
	$$ \left< e^{-\frac{t}{2} \Delta^1}(\omega_{\upsilon}) \cdot \eta_k \right> = \left< \omega_{\upsilon} \cdot e^{-\frac{t}{2} \Delta^1}( \eta_k )\right> \ . $$

Therefore, $\eta_k$ being an eigenform of eigenvalue $\lambda_k$, we have
 	$$e^{-\frac{t}{2} \Delta^1}( \eta_k ) = e^{- \frac{t}{2} \lambda_k} \eta_k \ , $$
 	
and thus
	$$ \left< e^{- \frac{t}{2} \Delta^1}(\omega_{\upsilon}) \cdot \eta_k \right> = e^{-\lambda_k \frac{t}{2} } \left< \omega_{\upsilon} \cdot \eta_k \right> \ . $$

The linking form $\Omega$ being integrable, one can use Fubini's theorem again to get 
	$$ \left< \omega_{\upsilon} \cdot \eta_k \right> = \inte{\upsilon} \left[ \inte{M} \ \Omega( (\cdot)^{\sharp} ,y)(\eta_k)_x \ d\mu_g(x) \right] \ .$$
	
But, by construction of the linking form as the operator $*dG$ kernel, we have 
\begin{align*}
	\inte{M} \ \Omega(x ,y)(\eta_k)_x \ d\mu_g(y) =
				 \left( (* d G) (\eta_k) \right)_y \ . 
\end{align*}	

The operator $G$ commutes with $d$, in particular all terms of this series corresponding to a closed differential form vanishes. The remaining terms are given by 
	$$ * d G(\eta_k) = \frac{*d\eta_k}{\lambda_k} \ .$$ 
	
Finally we have
		$$  \left< \omega_{\upsilon} \cdot \eta_k \right> =  \inte{\upsilon} *\left( \frac{d \eta_k}{\lambda_k} \right) \ ,$$
		
which gives 
	\begin{equation}
		 \left< e^{-\frac{t}{2} \Delta^1}(\omega_{\upsilon}) \cdot \eta_k \right>  = e^{-\frac{t}{2} \lambda_k} \inte{\upsilon} *\left( \frac{d \eta_k}{\lambda_k} \right) \ ,
	\end{equation}

the expected result. \hfill $\blacksquare$

\section{Application to torus geodesics linking.}

\label{4}

This section is devoted to the use of Formula \ref{theo formule localisebis} to compute the linking number of homologically trivial multi-geodesics of the canonical 3-tore torus $\mathbb{T}^3 := \quotient{\mathbb{R}^3}{\mathbb{Z}^3}$, for which the spectral theory is fully understood: we will describe it in Subsection \ref{soussec theorie spectrale}. One would like to use it to give a more or less explicit expression to the following series, $t > 0$ being fixed,
\begin{equation*}
	\label{equation series}
	 \somme{k >0} e^{- \lambda_k t} \inte{\gamma} \eta_k \inte{\upsilon} * \left( \frac{d \eta_k}{\lambda_k} \right) \ .
\end{equation*}

Right after, we will identify the limit of this series when $t \to 0$ to the Fourier development of some function. In the meantime, we will recall Theorem \ref{theo formule localisebis} which guarantees that this sequence of series actually converges to the linking number.

\subsection{Statement of the generalised torus linking theorem.}

\label{4.1}

We call a \textbf{mutli-geodesic} a multi-curve consisting of geodesics. This subsection goal is to state a formula giving the linking number of any two collections of multi-geodesics of the $\mathbb{T}^3$. \\

Let us fix some notations. We parametrise geodesics of $\mathbb{T}^3$ as follows 
\begin{equation*}
	\label{parametrisation}
	\left\{
	 \fonction{\gamma}{\mathbb{R}} { \mathbb{T}^3 }{t}{
	  \begin{pmatrix}
	 \gamma_1 t + \nu_1 \\ 
	 \gamma_2 t + \nu_2 \\
	 \gamma_3 t + \nu_3 
	 \end{pmatrix} \mod \mathbb{Z}^3 \ . }
	 \right. 
\end{equation*}

Note that these curves are automatically oriented by the parametrisation. \\

We call the point $\nu = \gamma(0) = \begin{pmatrix}\nu_1 \\ 
	 \nu_2 \\
	\nu_3
	 \end{pmatrix} \in \mathbb{T}^3 $ the \textbf{origin} of $\gamma$. The vector $\gamma'(t) \in T_{\gamma(t)} \mathbb{T}^3$ being also an element of the Lie algebra of $\mathbb{T}^3$ - canonically identified to $\mathbb{R}^3$ - one can write then $\gamma'(t)$ as
$$\gamma'(t)  = \begin{pmatrix}
	 \gamma_1  \\ 
	 \gamma_2  \\
	 \gamma_3  
	 \end{pmatrix} \in \mathbb{R}^3 \ .$$
	 
The vector $\gamma'(t)$ does not depend nor on $t$ neither on this origin of the curve. To each geodesic of $\mathbb{T}^3$ is therefore associated a vector of $\mathbb{R}^3$, which has the property that it belongs to $\mathbb{Z}^3$ if and only if $\gamma$ is closed. Note too that in this case any points $\nu \in \text{Im} (\gamma)$ can be chosen as origin. \newline
	 
\begin{figure}[!h]
	\begin{center}
	\label{figure torre}
		\def\svgwidth{0.4 \columnwidth}
		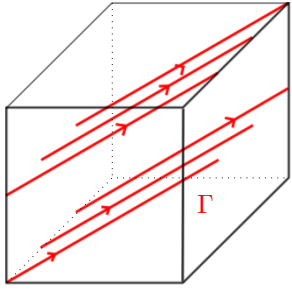
	\caption{The cube is identified to $\mathbb{T}^3$ by gluing opposed faces with translations. The red curve $\gamma$ admit as derivative vector $ \gamma' =  (2,2,3)$. This curve being closed, one can choose whatever of its points as an origin, for instance $0_{\mathbb{T}^3}$.}
		\end{center}
\end{figure}	 

The 3-torus fundamental group being Abelian, one can check that the vector $\gamma' \in \mathbb{Z}^3$ is canonically identified to the homology class of the closed curve $\gamma$ in $\mathbb{Z}^3$. We thus note $ [ \gamma ]$ the vector $\gamma'$ to emphasis its topological flavour. Off that remark comes out the following necessary and sufficient condition for a multi-geodesics $\Gamma = (\gamma^i)_{i \in I}$ to be  homologically trivial :
$$ \somme{i \in I} \left[ {\gamma^i} \right] = 0_{\mathbb{R}^3} \ .$$

The following construction is needed to state our theorem. \\

Given two vectors $u,v \in \mathbb{Z}^3$, we define the vector $\beta^{u, v} \in \mathbb{Z}^3$ as the unique one verifying the following conditions:

\begin{itemize}
\item $\beta^{u,v} \in \spn(u,v)^{\perp}$;
\item $\det(u, v, \beta^{u, v}) >0$;
\item its euclidean norm $|| \beta^{u, v} ||$ is minimal for the two first properties.
\end{itemize} 

Given two geodesics $\gamma$ and $\nu$, we still simply denote by $\beta^{\gamma, \nu}$ the vector $\beta^{[\gamma], [\nu]}$. Our torus linking theorem can then be stated as follows.

\begin{theoreme}
\label{thelktorus}
Let $\Gamma = (\gamma^i)_{i \in I}$ and $\Upsilon = (\upsilon^j)_{j \in J}$ two homologically trivial multi-geodesics of $\mathbb{T}^3$. They link according to the following formula: 

\begin{equation}
\label{fortheo lktorus}
\lk(\Gamma, \Upsilon) = \somme{i \in I , j \in J} \det \left( [\gamma^i], [\upsilon^j], \frac{\beta^{i,j}}{||\beta^{i,j}||}  \right) \frac{1 - 2 \lfloor (\nu^{i,j} \cdot \beta^{i,j}) \rfloor }{2 || \beta^{i,j} ||} 
\end{equation} \\

where $ \nu^{i,j} = {\gamma^i}(0) - \upsilon^j(0)$ is the difference between the two origins and $\lfloor \alpha \rfloor$ denotes the unique representative in $[0,1)$ of  the class $(\alpha  \mod \mathbb{Z})$. 
\end{theoreme}

\begin{remark} \
	\label{rem formule lk}
	\begin{itemize}
		\item One can define the linking number in every dimension $n$, providing that we considered two homologically disjointed sub-manifolds of dimension $p$ and $q$ satisfying $p+q = n-1$. Our method is likely to be generalised for flat torus in any dimension.
		\item \textit{A priori}, Formula \eqref{fortheo lktorus} depends on a choice of parametrisation. We will clarify this point along the proof with Remark \ref{rem parametrization}.
	\end{itemize}
\end{remark}

\subsection{Spectral theory of 1-differential forms of $\mathbb{T}^3$} 
\label{soussec theorie spectrale} 
\label{4.2}
We start by introducing some notations.

\begin{itemize} 
\item   We denote by a lower index $i$ the i-nth coordinate of a vector and by an upper index its belonging to a family of vectors. For example, $\gamma_i^j$ denotes the $i$-nth coordinate of the $j$-nth vector of a family indexed by $j \in J$. 
\item Given a vector 
	$$ v = \begin{pmatrix}
			v_1 \\ v_2 \\ v_3
			\end{pmatrix} 
		 \in \mathbb{R}^3 \ ,$$ 
	we note $v^*$ the differential form $ v_1 dx_1 + v_2 dx_2 + v_3 dx_3 $. This one being invariant by translations, it defines a harmonic differential form on $\mathbb{T}^3$. \textbf{We continue to note $v^*$ the induced-on-$\mathbb{T}^3$ differential form.}
\item The scalar product of two vectors $a$ and $b$ in $\mathbb{R}^3$ is denoted by $(a \cdot b)$ and the associated Euclidean norm by $||\cdot||$. 
\item The $\mathbb{R}^3$ vectorial product is denoted by $\wedge$. \\
\end{itemize} 

Let describe the 1-differential forms spectral theory of $\mathbb{T}^3$ thanks to the following set of datum:

\begin{itemize}
	\item a vector $ \textbf{k} = \begin{pmatrix}
k_1 \\ k_2 \\ k_3 
\end{pmatrix} \in \mathbb{Z}^3 $
	\item an orthonormal basis $(v^1, v^2, v^3)$ of $\mathbb{R}^3$.
	\item a function $f \in \{ \cos, \sin \}$
\end{itemize}   
 
\textit{ Note that here we have the choice on an orthonormal basis of $\mathbb{R}^3$.
}
\begin{lemma}
	\label{lemspectraltorus}
The 1-differential form of $\mathbb{T}^3$ 
			\begin{equation}
				\label{spectretore}
				\eta(x) = \sqrt{2} f ( 2 \pi (\textbf{k} \cdot x ) )(v^i)^* 
			\end{equation}
is an eigenform of $\Delta^1$ with associated eigenvalue $\lambda = (2 \pi ||\textbf{k}||)^2$. 
\end{lemma}

\textbf{Proof: } we start in showing that these forms are of unit norm: 
		\begin{align*}
	 		||\eta||_{L^2} :=& \inte{\mathbb{T}^3} \eta \wedge * \eta \\
	 		= & \inte{\mathbb{T}^3} 2 f ^2 ( 2 \pi (\textbf{k} \cdot x ) )(v^i)^* \wedge * (v^i)^* \\
	 		= & \inte{\mathbb{T}^3} 2 f ^2 ( 2 \pi (\textbf{k} \cdot x ) ) d \vol = 1 \ ,
		\end{align*}
since  $ f^2 = \frac{1 \pm f( 2 \cdot) }{2} $. \newline

Recall the Laplace operator definition $$ \Delta \eta =  (d \delta + \delta d) \eta \ .$$ 

Because $\delta = -*d* $ in dimension 3, one gets:
\begin{equation}
	d \delta \eta  = d(-*d*) \eta = - \sqrt{2} (d * d) \left( f ( 2 \pi (\textbf{k} \cdot x ) ) * (v^i)^* \right) .
\end{equation}

By the Hodge star definition we have  $$*(v^i)^* = (v^j)^* \wedge (v^t)^*$$

where $(i,j,t)$ is a circular permutation of $(1,2,3)$, so that 
\begin{equation}
			d \delta \eta  = - \sqrt{2} d * d \left( f ( 2 \pi (\textbf{k} \cdot x ) ) \wedge  (v^j)^* \wedge (v^t)^* \right) \ .
\end{equation}

And then,
\begin{align*}
		d \delta \eta & = - \sqrt{2} (d *) \left( 2 \pi k_i f'( 2 \pi (\textbf{k} \cdot x )  \left( (v^i)^* \wedge (v^j)^* \wedge (v^t)^* \right) \right)  \\
		& = -2 \sqrt{2}  \pi k_i d  \left( f'( 2 \pi (\textbf{k} \cdot x ) \right)  \\
		& = - 2 \sqrt{2}  \pi k_i  d f'( 2 \pi (\textbf{k} \cdot x ) \\
		& = - 4 \sqrt{2}  \pi^2  \left( k_i^2 f''( 2 \pi (\textbf{k} \cdot x ) dx_i + k_i k_j f''( 2 \pi (\textbf{k} \cdot x ) dx_j + k_i k_t f''( 2 \pi (\textbf{k} \cdot x ) dx_t \right) \ .
\end{align*}

We compute $ \delta d \eta $ in a similar way to get
\begin{align*}
\delta d \eta  = & - 4 \sqrt{2}  \pi^2  \left( k_j^2 f''( 2 \pi (\textbf{k} \cdot x ) dx_i + k_t^2 f''( 2 \pi (\textbf{k} \cdot x ) \right) \\
						&  + 4 \sqrt{2}  \pi^2 \left( k_i k_j f''( 2 \pi (\textbf{k} \cdot x ) dx_j - k_i k_t f''( 2 \pi (\textbf{k} \cdot x ) dx_t \right) \ ,
\end{align*} 
Summing both terms gives 
$$\Delta \eta = - 4 \sqrt{2}  \pi^2  \left( k_i^2 f''( 2 \pi (\textbf{k} \cdot x ) dx_i + k_j^2 f''( 2 \pi (\textbf{k} \cdot x ) dx_i + k_t^2 f''( 2 \pi (\textbf{k} \cdot x ) \right) dx_i \ .$$

Since $f'' = -f$ one has
$$\Delta  \eta = 4 \pi^2(k_1^2 +k_2^2 +k_3^2) \eta \ ,$$

the expected outcome. \hfill $\blacksquare$ \\

To use Theorem \ref{theo formule localisebis} we need a basis of eigenform. Fixing an orthonormal basis of $\mathbb{R}^3$, the family issued from all $\textbf{k} \in \mathbb{Z}^3$ and both the function $\cos$ and $\sin$ forms a generating family. To see it, one can decomposes in Fourier series the coefficients of a 1-differential form $\omega$ written as:
	$$ \omega(x,y,z) = f_1(x,y,z)  \cdot (v^1)^* + f_2(x,y,z)  \cdot (v^2)^* + f_3(x,y,z)  \cdot (v^3)^* \ .$$ 

Moreover, this family is free up to the trivial relations $\cos(-\textbf{k} \cdot x ) = \cos(\textbf{k} \cdot x )$ and $\sin(-\textbf{k} \sin x ) = -\sin(\textbf{k} \cdot x )$. 

\subsection{Computation of the approximating series} 
\label{subsec computation}

\label{4.3}

Recall that we parametrised both geodesics $\gamma$ and $\upsilon$ as:

\begin{equation*}
	\left\{ \fonction{\gamma}{ \quotient{\mathbb{R}}{\mathbb{Z}}}{ \mathbb{T}^3 }{t}{
	  	\begin{pmatrix}
			 \gamma_1 t + \nu_1 \\ 
			 \gamma_2 t + \nu_2 \\
	 		\gamma_3 t + \nu_3 
	 	\end{pmatrix}  }	\right. 
			 \hspace{0,5 cm}
	\left\{ \fonction{\upsilon}{ \quotient{\mathbb{R}}{\mathbb{Z}}}{ \mathbb{T}^3 }{t}{
	 	 \begin{pmatrix}
	 		\upsilon_1 t + \mu_1 \\ 
	 		\upsilon_2 t + \mu_2 \\
	 		\upsilon_3 t + \mu_3 
	 	\end{pmatrix}  } \right.
\end{equation*}

where $\gamma_i, \upsilon_j \in \mathbb{Z}$ and $\mu_j, \nu_j \in [0,1]$. \newline

First, note that we can assume that $\nu =0 $. In fact, since the translation of $\mathbb{R}^3$
 $$\tau_{\nu} := x \to x + \nu \ , $$ 
descends to an isometry of $\mathbb{T}^3$, using Corollary \ref{translation} one has:
	$$\inte{ \gamma} \inte{ \upsilon} \Omega = \inte{ \tau^{-1}(\gamma) } \inte{ \tau^{-1}(\upsilon)} \Omega \ ,$$

where now $  \begin{pmatrix} 0 \\ 0 \\ 0  \end{pmatrix}   $ belongs to $\tau^{-1}(\gamma)$. \textbf{In order not to burden the notations we will still denotes by $ \mu$ the new origin} (keeping in mind that it actually corresponds to $ \mu(\gamma,\upsilon) =\mu - \nu $) of the translated curve $\upsilon$. \newline

We saw that, given an orthonormal basis of $\mathbb{R}^3$, one can build an orthonormal eingenform basis of the Laplace operator. To simplify the calculation we will perform in \ref{equation localisebis} we make a choice of this orthonormal basis adapted to the curve $\gamma$: the first vector is chosen to be $v^1 = \frac{[\gamma]}{|| [\gamma] ||} $, and we arbitrary complete it to get an orthonormal basis:
			$$ \left(  \Big( v^1 = \frac{[\gamma]}{|| [\gamma] ||} \Big) ^*, v^2, v^3 \right)  \ .$$

Recall that we want to compute the following series;
\begin{equation}
	\label{equation series}
	 \somme{k >0} \ e^{- \lambda_k s} \inte{\gamma} \eta_k \inte{\upsilon} * \left( \frac{d \eta_k}{\lambda_k} \right) \ .
\end{equation}

We will compute all terms involved in this series separately and we will sum them in the next subsection. This terms are product of two integrals
%%\begin{equation}
%%\label{product}
%% \frac{e^{-\lambda t}}{\lambda} \left(\inte{\gamma} * d\eta \right) \left( \inte{\upsilon} \eta \right). 
%% \end{equation}
that we compute independently.\\
 
We start with the integral involving the operator $*d$. Let $\eta$ be an eingenform, one has 
\begin{equation}
	\label{equationpremier terme produit}
	 \inte{\gamma} \eta  = \inte{[0,1]} \sqrt{2} f( 2 \pi \textbf{k} \cdot \gamma(t) ) (v^i)^* ([\gamma]) dt \ ,
\end{equation} 

where $f \in \{ \cos, \sin \}$ and \textbf{k} is a vector of $\mathbb{Z}^3$. The above integral vanishes whenever $f$ is a sinus since the curve $\gamma$ passes by $0$;
\begin{align*}
	\inte{\gamma} \eta & = \inte{[0,1]} \sqrt{2} \sin( 2 \pi \textbf{k} \cdot \gamma(t) ) (v^i)^* ([\gamma]) dt \\
	& = C_1 \inte{[0,1]} \sin ( 2\pi C_2 t)dt	=0  
\end{align*}

because $C_2 \in \mathbb{Z}$. \textbf{We can assume then that $f$ is a cosinus.} We keep computing in considering the eigenforms
 		$$ 	\eta_{\textbf{k},i} = \sqrt{2} \cos ( 2 \pi (\textbf{k} \cdot x ) )(v^i)^* $$

only, where $ \textbf{k} \in \mathbb{Z}^3$ and $i \in \{1,2,3\}$. Which, looking backward to \eqref{equationpremier terme produit}, gives
$$ \inte{\gamma} \eta_{\textbf{k},i} = \inte{t = 0}^{1} \cos (2 \pi t (\textbf{k} \cdot [\gamma])) (v^i)^*( [\gamma] ) dt $$ 

where $(v^i)^*( [\gamma] ) =  ([\gamma] \cdot v^i) = ||[\gamma]|| \delta_{i,1}$. \newline

The above integral therefore vanishes whenever

\begin{itemize}
\item $(\textbf{k} \cdot [\gamma]) \neq 0$ ;
\item $ i \neq 1$ . 
\end{itemize} 

Moreover, in the case where it doesn't, the function $t \mapsto \cos \left((2 \pi (\textbf{k} \cdot [\gamma]) t \right) $ is constant, so that
\begin{equation}
\label{valeur1}
 \inte{\gamma} \eta_{\textbf{k},i} =  \sqrt{2} ||[\gamma]|| \ .
\end{equation}

Differential forms giving a non vanishing term of the series \eqref{equation series} are therefore 
$$\eta_{\textbf{k},1} = \sqrt{2} \cos( 2 \pi  \textbf{k} \cdot x ) \left( \frac{[\gamma]}{||[\gamma]||} \right)^* \ ,$$ 

with $\textbf{k}\in \mathbb{Z}^3$ and $\textbf{k} \cdot [\gamma]=0$.

\textbf{From now, we will denote $\eta_{\textbf{k},1}$ by $\eta_{\textbf{k}}$}. We now compute the second term of the series:
			$$ \inte{\upsilon} * d \eta_{\textbf{k}} \ ,$$ 
			
starting off computing $*d\eta_{\textbf{k}}$. 
\begin{align*}
	* d \eta_{\textbf{k}} &= *d \left(  \sqrt{2} \cos( 2 \pi  (x \cdot \textbf{k})) (v^1)^* \right) \\
	&=   -  2 \sqrt{2} \pi  \sin( 2 \pi (x \cdot \textbf{k}))* \left( k_1 dx_1 \wedge (v^1)^* + k_2 dx_2 \wedge (v^1)^* + k_3 dx_3 \wedge (v^1)^* \right) \\
	&=  -   2 \sqrt{2} \pi  \sin( 2 \pi  (x \cdot \textbf{k}))(\textbf{k} \wedge v^1)^*  \ .
\end{align*}

We then get 
\begin{align*}
\inte{\upsilon} * d\eta_{\textbf{k}} &= \int_{t=0}^1 - 2 \sqrt{2} \pi  \sin \left( 2 \pi t ([\upsilon] \cdot \textbf{k}) + 2 \pi ( \mu \cdot \textbf{k})  \right) (\textbf{k} \wedge v^1)^*([\upsilon]) dt \\
&=  - 2 \sqrt{2} \pi \det \left( \frac{[\gamma]}{||[\gamma]||}, [\upsilon], \textbf{k} \right) \int_{0}^1 \sin \left( 2 \pi t ([\upsilon] \cdot \textbf{k}) + 2 \pi ( \mu \cdot \textbf{k}) \right) dt \ . \\
\end{align*}

As well as before, this integral vanishes if one of this conditions holds: \\

\begin{itemize}
\item the vectors $[\gamma]$ and $[\upsilon]$ are collinear;
\item $(\textbf{k} \cdot [\upsilon]) \neq 0$. \\ 
\end{itemize}

Moreover if $\inte{\nu}  * d \eta_{\textbf{k}} \neq 0$ we have
\begin{equation}
\label{valeur2}
 \inte{\upsilon} * d \eta_{\textbf{k}} = -  2 \sqrt{2} \pi \det \left( \frac{[\gamma]}{||[\gamma]||}, [\upsilon], \textbf{k} \right) \sin(2 \pi ( \mu \cdot \textbf{k})) \ .
 \end{equation}

Multiplying \eqref{valeur1} and \eqref{valeur2} one has: \\
\begin{equation*}
	\inte{\gamma} \eta_{\textbf{k}} \inte{\upsilon} * d\eta_{\textbf{k}} = 
		\left\{ 
			\begin{array}{l}
				4 \pi \det( [\gamma], [\upsilon], \textbf{k})  \sin(2 \pi ( \mu \cdot \textbf{k})) \ \text{ si } \ \textbf{k} \in \text{Span}([\gamma], [\upsilon])^{\perp}  \\
				0 \ \text{ sinon} 
			\end{array} \right.
\end{equation*}

This leads us to characterise elements of $\text{Span}([\gamma], [\upsilon])^{\perp} $ with integer coefficients.

\begin{lemma} Let $b_1$ and $b_2$ two non zero vectors of $\mathbb{Z}^3$ then the group
	\begin{equation*}
 		\mathrm{Span}(b_1, b_2)^{\perp} \cap \mathbb{Z}^3 
	\end{equation*}
is cyclic. We note by $\pm \beta$ one of this two possible generators.
\end{lemma}

\textbf{Proof :} As a set it is non empty; the vector $b_1 \wedge b_2$ belongs to $\mathbb{Z}^3$ and is orthogonal to both $b_1$ and $b_2$. As the intersection of two subgroups, $\mathbb{Z}^3$ and $\mathbb{R} \cdot b_1 \wedge b_2$, it is a subgroup of $\mathbb{R}$. The neutral element of  		$\mathrm{Span}(b_1, b_2)^{\perp} \cap \mathbb{Z}^3$ must be isolated since $\mathbb{Z}^3$ is discrete. By characterisation of $\mathbb{R}$ subgroups, this group is cyclic. \hfill $\blacksquare$ \\

We apply the previous lemma to the pair $([\gamma], [\upsilon])$ to get the following description of elements $\textbf{k} \in \mathbb{Z}^3$ giving a non vanishing term in the series of Equation \ref{equation localisebis}:
\begin{equation*}
 \text{Span}([\gamma], [\upsilon])^{\perp} \cap \mathbb{Z}^3 = \{ k \ \beta, \hspace{0.2cm} k \in \mathbb{Z} \}  \ .
\end{equation*}

Among both possible generators, we choose $\beta$ such that the family $([\gamma], [\upsilon], \beta)$ is positively oriented. \\

The only non vanishing terms of the series appearing in \eqref{equation series} correspond to the differential forms 
	$$ \eta_{ (k \beta) } = \sqrt{2} \cos \big( (k \beta) \cdot x \big) \left(\frac{[\gamma]}{||\gamma||} \right)^* \ , $$

and in this case we have
\begin{equation*}
 	\inte{\gamma} \eta_k \inte{\upsilon} * d\eta_k = - 4 \pi k \det( [\gamma], [\upsilon], \beta)  \sin(2 \pi ( k (\mu \cdot \beta))) \ .
\end{equation*}

\textbf{From now we denote by $\eta_{k}$ the differential form $\eta_{ (k \beta) }$.} Recall that differential forms $\eta_k$ and $\eta_{-k}$ are collinear. To get a free family of eingenforms one needs to choose the sign of the integers $k$: we take them non negative. The series of Equation \ref{equation localisebis} then becomes: \\
\begin{equation}
	\label{equation ref grosse serie}
  - \somme{k >0} \frac{e^{- (2 \pi ||\beta||)^2 n^2 s}}{ \pi k ||\beta||  } \det \left( [\gamma], [\upsilon], \frac{\beta}{||\beta||} \right)  \sin(2 \pi  k (\mu \cdot \beta))) \ . 
\end{equation} \\

\begin{remark}
\label{rem parametrization}
As noticed in \ref{rem formule lk}, Formula \ref{equation ref grosse serie} is not  independent of the parametrisations involved \textit{a priori}. In fact, the point $\mu \in \mathbb{T}^3$ appearing in the term $\sin(2 \pi  k (\mu \cdot \beta)))$ depends of an origin choice for $\upsilon$. Let us thus check that $k (\mu \cdot \beta))$ actually doesn't, modulo $\mathbb{Z}$. Let $\mu_2 \in \upsilon $ an other origin of $\upsilon$, by definition there is $ t \in \mathbb{R}$ and $\alpha \in \mathbb{Z}^3$ such that $$ \mu_2 - \mu = t [\upsilon] + \alpha \ ,$$ 
thus
$$(\mu_2 \cdot \beta) = (\mu_2 - \mu + \mu \cdot \beta) = (\mu \cdot \beta) + (\alpha \cdot \beta) \ ,$$
since $ \beta \in [\upsilon]^{\perp} $. We conclude reducing the above formula modulo $\mathbb{Z}$ to get 
	$$(\mu_2 \cdot \beta) = (\mu \cdot \beta) \ , $$ 
since $(\alpha \cdot \beta) \in \mathbb{Z}$.
\end{remark}

 \subsection{A uniformly converging family of functions.}
\label{4.4}
Let us now look into the series \eqref{equation ref grosse serie} more in detail. If one is able to let $t \to 0$ within all terms of this series one would get
$$ - C \ \somme{k >0} \frac{1}{k} \sin(2 \pi k x) \ ,$$

with $C =   \frac{1}{ \pi ||\beta||}  \det \left( [\gamma], [\upsilon], \frac{\beta}{||\beta||} \right) $ and $ x = (\mu \cdot \beta)$. \newline

On can recognise here the Fourier series development of the defined-on-the-circle-$\quotient{\mathbb{R}}{\mathbb{Z}}$ function
\begin{equation}
	x \mapsto 
		\left\{ 
			\begin{array}{l l }
				0 &\text{ if } x = 0 \\
				\frac{\pi}{2}( 1 - 2x ) & \text{ on } (0, 1)
			\end{array}
		\right.
\end{equation}

So that we would have 
\begin{equation}
\label{formule}
\inte{\gamma} \inte{\upsilon} \Omega =  \frac{1}{2||\beta||}  \det \left( [\gamma], [\upsilon], \frac{\beta}{||\beta||} \right) (1 - 2(\mu \cdot \beta)) \ ,
\end{equation} 

which is precisely what expected. To justify the term-by-term convergence of \eqref{equation ref grosse serie}, we use the following lemma 

	\begin{lemma}
		\label{lemme series de fonction}
		Let $a_k(t) $ and $b_k(t)$ two sequences of functions defined on an interval $I$ containing $0$ such that
		\begin{enumerate}
			\item  $\left( \somme{ k \le n} a_k(t)\right)_{n \in \mathbb{N}}$ is uniformly bounded with respect to $t$;
			\item  the sequence of function $ b_k(t)$ is non-increasing with respect to $t$ and converges uniformly ,with respect to $k$, to $0$,
		\end{enumerate}		 
	then the series of functions $\somme{k \in \mathbb{N}} a_k(t) b_k(t)$ converges uniformly on $I$.
	\end{lemma}

We omit the proof, which consists to a discrete integration by part of the series. \\

We set $a_k(t) = \sin(2 \pi k x)$, $b_k(t) = \frac{e^{- a  t k^2}}{k}$ and $I = [0, + \infty]$. One can then check that for all $a \in \mathbb{R^+}$ and $x>0$, all assumptions of Lemma \ref{lemme series de fonction} holds. We deduce that the series of functions
$$ \somme{k >0} \frac{ e^{-a t k^2 } }{k}  \sin(2 \pi k x) $$

converges uniformly on $]0,+ \infty]$. One is therefore allowed to switch limits and sum in Equation \ref{equation ref grosse serie} to get:
\begin{align*}
	\limi{t \to 0} \ \somme{k >0} \frac{ e^{-a t k^2 } }{k}  \sin(2 \pi k x) & =  \somme{k >0} \ \limi{t \to 0} \  \frac{ e^{-a t k^2 } }{k}  \sin(2 \pi k x) \\
	& = \somme{k >0} \frac{\sin(2 \pi k x)}{k} \ ,
\end{align*}

which concludes the proof. 
\subsection{The $\mathbb{T}^2$-geodesic flow special case}

\label{4.5}

Particularly interesting collections of multi-geodesics of $\mathbb{T}^3$ arise as periodic orbits of the $\mathbb{T}^2$-geodesic flow. More generally, linking number of collections of periodic orbits have been studied by E. Ghys \cite{articleghysknots} and P. Dehornoy \cite{articlepierrelefthanded} \cite{thesedehornoy} for dynamical purposes. The later showed that, for al large class of examples given by geodesic flows on surfaces, these collections all link positively, up to a choice of global orientation. This implies the existence of Birkhoff sections and, as a corollary, that periodic orbits of this flows display fibred knots. In the setting of $\mathbb{T}^2$, we will see that Theorem \ref{thelktorus} specifies easily giving a new linking number formula. \\

We start by noticing that $\mathbb{T}^3$ is identified to the unitary tangent bundle $U \mathbb{T}^2$ of the 2-torus $\mathbb{T}^2 := \quotient{\mathbb{R}^2}{\mathbb{Z}^2}$. In fact, the unitary tangent bundle of $\mathbb{T}^2$ is trivial, $\mathbb{T}^2$ being a Lie group. One trivialisation consists to choose a direction of $\mathbb{R}^2$, which induces one on  $\mathbb{T}^2$, from which one is able to assign an angle to any vector of $U \mathbb{T}^2$. That is to say the map
 	$$ \fonctionbis{ U \mathbb{T}^2 }{\mathbb{T}^3}{ u}{ \left( (x,y), \theta \right) }  $$
is an actual trivialisation. With the unitary tangent bundle of a Riemannian manifold comes always a flow: the geodesic flow. In the case of $\mathbb{T}^2$, one can fully write down the flow in the trivialisation given above  
\begin{equation*}
	\label{fonction}
		\left\{ \fonction{\Phi_t}{ \mathbb{T}^3}{ \mathbb{T}^3}{ (x,y, \theta) }{ (x + t \cos \theta ,y + t \sin \theta, \theta) \ .} \right. 
\end{equation*}

Note that periodic orbits of a flow are naturally parametrised and oriented by the flow itself:
\begin{equation*}
		\left\{ \fonction{\gamma}{ \mathbb{S}^1}{ \mathbb{T}^3}{t}{ (x + t \cos \theta ,y + t \sin \theta, \theta) \ . } \right. 
\end{equation*}
 
\begin{figure}[!h]
\centering
\includegraphics[scale=0.35]{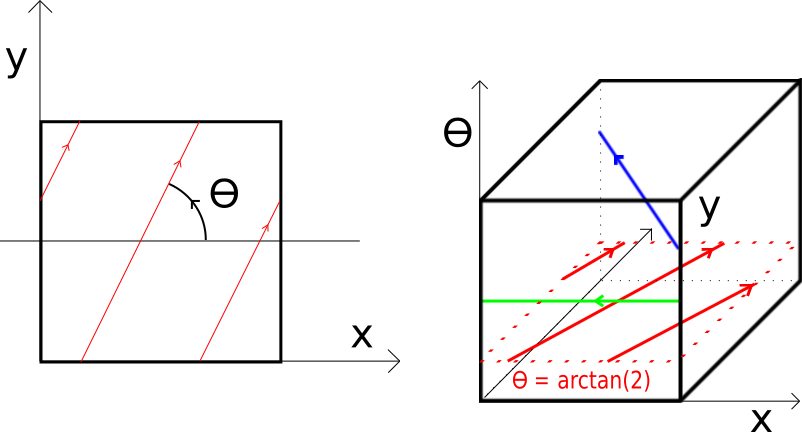}
\caption{The red curves on the left represents a closed geodesic of $\mathbb{T}^2$. This curve lifts canonically to the the red right one on the unitary tangent bundle. This lifted curve remains in the leaf $\theta = \arctan (2)$. Both blue and green curves represent two others lifted geodesics.}
\label{figure geo}
\end{figure}

\begin{remark}
The fact that orbits of the geodesic flow are still geodesics on the unitary tangent bundle is more general, providing that one endows the later with the right metric; the so called Sasaki metric. In our case, it turns out that the Sasaki metric coincides with the $\mathbb{T}^3$ flat one.
\end{remark}

In this setting, one can readily specifies Theorem \ref{theo formule localisebis} to get the

\begin{corollary}\cite[page 11]{thesedehornoy}
\label{corollaire Pierre2}
Let $\Gamma = \gamma^i$ and $\Upsilon = \upsilon^j$ two homologically trivial multi-geodesics of $\mathbb{T}^2$. In the $\mathbb{T}^2$-unitary tangent bundle they link according to the following formula:
$$ \lk (\Gamma, \Upsilon) = \somme{i \in I,j \in J}   i ( \gamma^i, \upsilon^j ) \frac{1 -  \frac{x_{i,j}}{\pi}}{2} \ ,$$
where $x_{i,j}$ denotes the unique determination in $[0,2\pi[$ of the oriented angle $\theta$ made at any intersections points (see Figure \ref{fig tore pierre}), and $i ( \gamma^i, \upsilon^j )$ denotes the algebraic intersections between $ \gamma^i$ and $\upsilon^j$ on $\mathbb{T}^2$.
\end{corollary}

\textbf{Proof :} as previously noticed, the orbits of this flow remain in the leafs $\theta = \text{cst} $, so that  the vectors $[{\gamma^i}]$ and $[{\upsilon^j}]$ belong $\mathbb{R}^2 \subset \mathbb{R}^3$. Our vector $\beta^{i,j}$ defined in Theorem \ref{thelktorus} becomes 
$$ \beta^{i,j} = \begin{pmatrix}	
0 \\ 0 \\ \pm 1
\end{pmatrix}$$ 
for all pairs $(i,j)$, the sign depending whether or not the angle between the curves $\gamma_i$ and $\upsilon_j$ is greater than $\pi$. In particular we have $ || \beta^{i,j} || =1$. Moreover, $\det( [{\gamma^i}], [{\upsilon^j}], \beta^{i,j})$ becomes $\det_{\mathbb{R}^2}( [\gamma_i], [\upsilon_j])$, which  corresponds to the intersection number between $\gamma^i$ and $\upsilon^j$ seen as curve of $\mathbb{T}^2$. To conclude, the quantity $ ( \beta^{i,j} \cdot \mu^{i,j} ) $ turns out to be interpreted as the difference between the angle made by the curve, \textit{i.e}: 
	$$(\pi - (\mu^{i,j} \cdot \beta^{i,j})) =  (\pi - (x_{i,j})) \ . $$
 \hfill $\blacksquare$

\bibliographystyle{alpha}

\bibliography{bibliography} 

\end{document}